\documentclass[reqno]{amsart}
\usepackage{cite}
\usepackage{mathrsfs}
\usepackage{color}
\usepackage{amsmath}
\usepackage{amsfonts}
\usepackage{amssymb}
\usepackage{graphicx}
\usepackage{hyperref}
\usepackage{multirow}
\usepackage{caption}
\usepackage{float}
\usepackage[all,pdf]{xy}
\usepackage{tikz}
\usepackage{tikz-cd}
\usepackage[numbers,sort&compress]{natbib} 
\usepackage[left=4.5cm,right=4.5cm,bottom=2.5cm]{geometry}

\newtheorem{Theorem}{Theorem}[section]

\newtheorem{Lemma}[Theorem]{Lemma}

\newtheorem{Definition}[Theorem]{Definition}

\newtheorem{Example}[Theorem]{Example}

\numberwithin{equation}{section}

\title{Codimension 1 transfer maps of K theoretic indexes}
\author{Yuetong Luo}
\begin{document}
	\setlength{\baselineskip}{15pt}
	\vspace*{-1cm}
	\begin{abstract}
		Let $M$ be a closed spin manifold and $N$ be a codimension 1 submanifold of it. In \cite{Zeidlercodim1}, given certain homotopy conditions, Zeidler shows that the Rosenberg index of $N$ is an obstruction to the existence of positive scalar curvature on $M$. He further gives a transfer map between the K groups of the group $C^*$ algebras of the foundemental group. The transfer map maps the Rosenberg index of $M$ to the one of $N$. In this note, we present an alternative formulation of the transfer map using maps between $C^*$ algebras, and give an analogus result for the codimension 1 transfer of higher K theoretic signatures.
	\end{abstract}
	\maketitle
	\begingroup
	\let\clearpage\relax
	\section{Introduction}
An important question in geometric topology is whether a given smooth manifold admits a Riemannian metric with positive scalar curvature. In the case of spin manifolds, the Sch{\"o}dinger-Lichnerowicz formula \cite{Schrodinger,Lichnerowiczformula} implies that the index of the Dirac operator vanishes on a compact spin manifold with positive scalar curvature, providing an obstruction to the existence of positive scalar curvature.

Rosenberg \cite{Rosenbergindex1,Rosenbergindex2,Rosenbergindex3} used the index of Dirac operators twisted with certain Hilbert module over $C^*$-algebras to obtain strong obstructions to the existence of positive scalar curvature. In particular, using the Mishchenko bundle on a spin manifold $M$, one can obtain the Rosenberg index $\alpha(M) \in K_*(C^*\pi_1(M))$.

Gromov and Lawson \cite{GromovLawsonpsc} found an intriguing obstruction to positive scalar curvature from submanifolds of codimension two with certain conditions. Hanke-Pape-
Schick \cite{HPScodim2} illuminate the result in the K-theoretic index way and Zeidler \cite{Zeidlercodim1} treats the case of codimension one submanifolds:

\begin{Theorem}[Theorem 1.7 in \cite{Zeidlercodim1}]
	\label{transfer}
	\
	
	Let $M$ be a closed connected spin manifold of dimension $m$ and $N \subset M$ a connected submanifold of codimension one with trivial normal bundle. Assume that the inclusion map $N \hookrightarrow M$ induces an injection on the fundamental group. 
	
	Then there is a homomorphism  (called the transfer map) for both reduced and maximal group $C^*$-algebras: $$\rho_{M,N}: K_*(C^*\pi_1(M)) \longrightarrow K_{*-1}(C^*\pi_1(N))$$
	s.t. $\rho_{M,N}$ maps the Rosenberg index of $M$ to the one of $N$, i.e. $\rho_{M,N}(\alpha(M))=\alpha(N)$.

\end{Theorem}

The main result of this paper is that the above transfer map can be realized by a certain homomorphism between $C^*$-algebras:

\begin{Theorem}
	\label{main thm}
	\
	
	Let $M$ with submanifold $N$ be as in Theorem \ref{transfer}. Let $H_{C^*\pi_1(N)} \stackrel{def}{=} l^2(\mathbb{N}) \otimes C^*\pi_1(N)$ be the standard
	Hilbert $C^*\pi_1(N)$ module, denote $K_{C^*\pi_1(N)}$ and $B_{C^*\pi_1(N)}$ to be its $C^*$-algebras of compact and bounded adjointable operators. Denote $\mathcal{Q}_{C^*\pi_1(N)}$
	$\stackrel{def}{=}B_{C^*\pi_1(N)}/K_{C^*\pi_1(N)}$ to be the Calkin algebra associated to $C^*\pi_1(N)$. 
	
	Then there is a homomorphism $\rho: C^*\pi_1(M) \longrightarrow \mathcal{Q}_{C^*\pi_1(N)}$ for both reduced and maximal group $C^*$-algebras, s.t. if we denote $\rho_{M,N}$ to the map $\rho_* \colon K_*(C^*\pi_1(M)) \longrightarrow K_*(\mathcal{Q}_{C^*\pi_1(N)}) \cong K_{*-1}(C^*\pi_1(N))$, then:
	
	(1) If $M$ is spin, then $\rho_{M,N}$ maps the Rosenberg index of $M$ to be the one of $N$, i.e. $\rho_{M,N}(\alpha(M))=\alpha(N)$.
	
	(2) Let $\epsilon=\left\{
	\begin{aligned}
		& 1  & if \ m \equiv 0 \mod 2\\
		& 0 & if \ m \equiv 1 \mod 2
	\end{aligned}
	\right.$. If $M$ is oriented, then $\rho_{M,N}$ maps the $C^*$-algebraic higher signature of $M$ to $2^{\epsilon}$ times the one of $N$. More generally, let $f:M' \longrightarrow M$ be a map between closed oriented manifolds of degree 1. Assume $N$ is transversal to $f$ and set $N'=f^{-1}(N)$, then
	$$\rho_{M,N}(sgn(M',f^*\nu_M))=2^{\epsilon}sgn(N',f^*\nu_N)$$

\end{Theorem}

Here is a brief outline of the paper. In Section 2 we shall introduce the general geometric setup that we work on. In Section 3 we shall construct the transfer map stated in Theorem \ref{main thm} in the maximal $C^*$-algebra case. In Section 4 we will prove that the transfer map constructed satifies the properties in Theorem \ref{main thm}. The proof is similar to the proof in \cite{HPScodim2}, with the constructions slightly different because of different geometric setups. In the last section we will extend the transfer map to the reduced $C^*$-algebra case.

The author would like to thank his advisor Thomas Schick for suggesting the research question and for carefully reading the first draft. The author is supported by DAAD Graduate School Scholarship Programme (57650678).

	\section{Geometric setting}
In this section, we will introduce the general geometric setup that will be used in the construction of the transfer map.

Geometric setup:

(1) $M$ is a closed and connected $m$ dimensional manifold, $N \subset M$ is a codimensional 1 connected submanifold with trivial normal bundle.

(2) Fix a base point $x_0 \in N \subset M$. Let $\pi_1(M,x_0)=\Gamma,\pi_1(N,x_0)=\pi$ and suppose that the inclusion map induces an injection $\pi \hookrightarrow \Gamma$ on the fundamental group.

(3) Denote $\tilde{p}:\widetilde{M} \longrightarrow M$ to be the universal cover of $M$ and let $\overline{M}=\pi \backslash \widetilde{M}$. Denote $\underline{p}:\widetilde{M} \longrightarrow \overline{M}$, $p:\overline{M} \longrightarrow M$ to be the corresponding covering maps. Fix $\tilde{x}_0 \in \widetilde{M}$ and set $\bar{x}_0=\underline{p}(\tilde{x}_0)$.

(4) Choose a tubular neighborhood of $N$: $t:N \times \mathbb{R} \longrightarrow M$ and lift it to the unique embedding $\bar{t}:N \times \mathbb{R} \longrightarrow \overline{M}$ with $\bar{t}(x_0)=\bar{x}_0$. Then $
\overline{M} \backslash \bar{t}(N \times \{0\})$ consists of two connected components, the closures
of these components are denoted by $\overline{M}_{\pm}$, where the sign in the notation is determined by
the requirement $\bar{t}(N \times \{\pm 1\}) \subset \overline{M}_{\pm}$. Denote $\widetilde{N}=\underline{p}^{-1}(\bar{t}(N \times \{0\}))$, $\widetilde{M}_{\pm}=\underline{p}^{-1}(\overline{M}_{\pm})$, they are path connected.

(5) The restriction of the universal cover map $\tilde{p}$ to $\widetilde{N}$: $\tilde{p}|_{\widetilde{N}}: \widetilde{N} \longrightarrow N$, is a universal cover for $N$. Lift $\bar{t}$ to the unique embedding $\tilde{t}: \widetilde{N} \times \mathbb{R} \longrightarrow \widetilde{M}$ with $\tilde{t}(x_0)=\tilde{x}_0$.
		\section{Transfer map}
	\label{transfer map}
	\label{constructiontransfer}
	In this section, we will construct the tranfer map stated in theorem \ref{main thm}. In the following sections till the last section, all group $C^*$–algebras
	are maximal group $C^*$–algebras. For simplicity, we will write $C^*G$ for the maximal group $C^*$–algebra of a group $G$.
	
	Here is a brief idea of how the construction goes, we consider the $C^*\pi$ module $\Gamma \times_{\pi} C^*\pi=\mathop{\oplus}\limits_{s \in \Gamma/\pi}C^*\pi$, $\Gamma$ acts left on it and we can complete the module to a Hilbert $C^*\pi$ module in a natural way. From geometric view $\Gamma/\pi$ can be identified with the set $p^{-1}(x_0) \subset \overline{M}$, thus we can consider how $\Gamma$ acts restricted to the parts corresponding to $\overline{M}_{+}$ and it turns out that in this way we can construct a homomorphism $\rho:\Gamma \longrightarrow \mathcal{Q}_{C^*\pi}$. The homomorphism extends to $C^*\Gamma$ by the universal property of the maximal group $C^*$-algebra.
	
	We begin with the construction. Fix an odd function $\chi_0:\mathbb{R} \longrightarrow [-1, 1]$ with $\chi_0(\pm s) = \pm 1$ for $s > 1$. We define a function $\chi: \overline{M} \longrightarrow [-1, 1]$ as follows:
	
\begin{equation*}
	\chi(\bar{x})=\left\{
	\begin{aligned}
		& \chi_0(s)  &\bar{x}&=\bar{t}(x,s) \in \bar{t}(N \times [-1,1])\\
		& \pm 1 &\bar{x}& \in \overline{M}_{\pm} \backslash \bar{t}(N \times [-1,1])
	\end{aligned}
	\right.
\end{equation*}

Fix $\tilde{x}_0 \in \widetilde{M}$, we define a function $\chi_{\Gamma / \pi}: \Gamma /\pi \longrightarrow [-1,1]$ as follows: for $\gamma \pi \in \Gamma/\pi$, we have $\chi_{\Gamma / \pi}(\gamma \pi)=\chi(\underline{p}(\gamma^{-1} \tilde{x}_0))$.

In order to give the description of the transfer map, we recall some basic definitions first:

\begin{Definition}[Hilbert module]
	\
	
	Let $A$ be any $C^*$-algebra, a Hilbert $A$ module $E$ is a complex vector space with right $A$ action and a sesquilinear map (linear on the second varible) $<-,->: E \times E \longrightarrow A$, called the $A$ valued inner product, such that for all $x,y \in E, a \in A$ we have:
	
	(1) $<y,x>=<x,y>^*$
	
	(2) $<x,x> \geq 0$ and $<x,x>=0$ if and only if $x=0$.
	
	(3) $<x,ya>=<x,y>a$
	
	(4) Denote $\lVert \cdot \rVert_A$ to be the $C^*$ norm on $A$, $E$ is complete with respect to the norm $\lVert x \rVert=\sqrt{\lVert <x,x> \rVert_A}$.
\end{Definition}

\begin{Example}[Standard Hilbert module]
	\
	
	Let $A$ be a $C^*$-algebra, denote $A^{\infty}$ to be the direct sum of countable many copies of $A$. The elements of $A^{\infty}$ can be written as $(a_1,a_2,...)$ with $a_i \in A$ and $a_i \neq 0$ for all but finite many $i$. There's a sesequilinear map on $A^{\infty}$ given by $<x,y>=\sum\limits_{i \geq 1}a_ib_i$ for $x=(a_1,a_2,..)$ and $y=(b_1,b_2,...)$. The sesequilinear map satisfies the property (1) to (3) in the definition. Thus if we complete the space with respect to the norm $\lVert x \rVert=\sqrt{\lVert <x,x> \rVert_A}$, then we get a Hilbert $A$ module. It is often called the standard Hilbert $A$ module and we will denote it by $H_A$.
\end{Example}

\begin{Definition}[bounded adjointable operators]
	\
	
	Let $A$ be a $C^*$-algebra and E be a Hilbert $A$ module. $B_{A}(E)$ is the set of all module
	homomorphisms $T:E \longrightarrow E$ such that there is an adjoint module homomorphism $T^*:E \longrightarrow E$ with $<Tx,y>=<x,T^*y>$ for all $x,y \in E$.
	
	We denote $B_A(H_A)$ as $B_A$ and $GL_A(E)$ to be the set of all invertible elements in $B_A(E)$.
\end{Definition}

Consider the vector space $\Gamma \times_{\pi} C^*\pi=\mathop{\oplus}\limits_{s \in \Gamma/\pi} C^*\pi$, there is a right action of $C^*\pi$ given by right multiplication and there is a $C^*\pi$ valued inner product given by $<a,b>=\sum\limits_{s \in \Gamma/\pi}a_s^*b_s$ for $a=(a_s) \in \mathop{\oplus}\limits_{s \in \Gamma/\pi} C^*\pi$ and $b=(b_s) \in \mathop{\oplus}\limits_{s \in \Gamma/\pi} C^*\pi$. Complete $\mathop{\oplus}\limits_{s \in \Gamma/\pi} C^*\pi$ with respect to the norm $||a||=\sqrt{||<a,a>||_{C^*\pi}}$ and denote the resulting Hilbert $C^*\pi$ module to be $E$. $\Gamma$ acts unitarily on $\Gamma \times_{\pi} C^*\pi$ by the induction of left action of $\pi$ on $C^*\pi$, thus the action extends continously to $E$. We denote the action map to be $g \in \Gamma \longmapsto \varphi_g \in B_{C^*\pi}(E)$.

For any $S \subset \Gamma / \pi$, we define $E_S$ as the subset of $E$ consisting of those elements in which $x_s=0$ for all $s$ not in $S$. Similarly, for any subset $W$ of $[-1,1]$, we define $E_W$ as $E_{\chi^{-1}_{\Gamma / \pi}(W)}$. The $E_S,E_{S^c}$ and $E_W,E_{W^c}$ are complemented submodules of $E$ respectively. We denote $I_W$ and $I_S$ to be the inclusion maps from $E_W$ and $E_S$ to $E$ respectively, and denote $P_W$ and $P_S$ to be the projection maps from $E$ to $E_W$ and $E_S$ respectively. 

Fix any number $r \in (-1,1)$, let $F_r=E_{[r,1]} \oplus H_{C^*\pi}$, where $H_{C^*\pi}$ is the standard Hilbert  $C^*\pi$ module. We define the map $\rho:\Gamma \longrightarrow \mathcal{Q}_{C^*\pi}(F_r)=\mathcal{Q}_{C^*\pi}$ by $g \longmapsto \rho_g=\varPi \circ (P_{[r,1]}\varphi_gI_{[r,1]} \oplus Id)$, where $\mathcal{Q}_{C^*\pi}(F_r)=B_{C^*\pi}(F_r)/K_{C^*\pi}(F_r)$ is the Calkin algebra of the Hilbert module $F_r$ over $C^*\pi$ and $\varPi$ is the natural projection map. Then it is easy to check that $\rho$ is a group homomorphism and $\rho_g^*=\rho_{g^{-1}}$. By the universal property of the
maximal group $C^*$-algebra, the homomorphism $\rho$ extends to a * homomorphism $\rho: C^*\Gamma \longrightarrow \mathcal{Q}_{C^*\pi}$.

We define the transfer map $\rho_{M,N}$ to be $\rho_*: K_*(C^*\Gamma) \longrightarrow K_*(\mathcal{Q}_{C^*\pi}) \cong K_{*-1}(C^*\pi)$. One basic property of the transfer map is that it's independent of the choice of $r$:

\begin{Lemma}
	\label{indelem}
	The map $\rho_{M,N}$ is independent of the choice of $r$.
\end{Lemma}

\begin{proof}
	The proof is given in Section \ref{pfindelem}.
\end{proof}

	\section{Proof of theorem \ref{main thm}}
We will prove theorem \ref{main thm} in this section. Due to the difference in describing the boundary map of Mayer-Vietoris sequences, the proof is divided into two cases:
$N$ is separating in $M$ and $N$ is non-separating in $M$.

Before getting to the main proof, we make the following observations:

(1) There is a Mayer-Vietoris exact sequence in K-homology for the decomposition: $M=N \times I \cup M_N$, where $N \times I=t(N \times [-2,2])$ and $M_N=M \backslash t(N \times (-1,1))$:
\begin{equation*}
	\begin{tikzcd}
		... \rar & K_*(N \times [-2,-1] \amalg N \times [1,2]) \rar \dar[phantom, ""{coordinate, name=Z}] & K_*(N \times I) \oplus K_*(M_N) \arrow[dll,
		rounded corners,
		to path={ -- ([xshift=2ex]\tikztostart.east)
			|- (Z) [near end]\tikztonodes
			-| ([xshift=-2ex]\tikztotarget.west)
			-- (\tikztotarget)}]\\
		K_*(M) \rar{\partial} & K_{*-1}(N \times [-2,-1] \amalg N \times [1,2]) \rar & ...\\
	\end{tikzcd}
\end{equation*}
Since $K_*(N \times [-2,-1] \amalg N \times [1,2])=K_*(N) \oplus K_*(N)$ and $K_*(N \times I)=K_*(N)$, by the exactness of the sequence, $Im\partial \subset \{(x,-x)|\ x \in K_{*-1}(N)\}$. Thus the $\partial$ map in the sequence gives a map $\partial_{MV}:K_*(M) \longrightarrow K_{*-1}(N)$. A spin structure on $M$ defines a fundamental K–homology class $[M] \in K_m(M)$ and standard result shows that $\partial_{MV}[M]=[N]$. If $M$ is oriented, there is a signature class on $M$: $[M]_{sgn} \in K_m(M)$ and $\partial_{MV}[M]_{sgn}=2^{\epsilon}[N]_{sgn}$ (see Proposition 4.1 in \cite{boundarysignature}), where

\begin{equation*}
	\epsilon=\left\{
	\begin{array}{rcl}
		1 & \text{if } m \equiv 0 \mod 2\\
		0 & \text{if } m \equiv 1 \mod 2
	\end{array}
	\right.
\end{equation*}

(2) The Mishchenko bundle $\nu_M=\widetilde{M} \times_{\Gamma} C^*\Gamma$ defines a class $[\nu_M] \in K^0(M;C^*\Gamma)=K_0(C(M) \otimes C^*\Gamma)$. The Rosenberg index $\alpha(M)$ is the pairing of the class $[M]$ with this class and the $C^*$-algebraic higher signature $sgn(M,\nu_M)$ is the pairing of the signature class with this class. 

(3) There is a short exact sequence involving the Calkin algebra:
$$
\begin{tikzcd}
	0 \rar & K_{C^*\pi} \rar & B_{C^*\pi} \rar& \mathcal{Q}_{C^*\pi} \rar & 0\\
\end{tikzcd}
$$ and its associated long exact sequence in K–theory yields the boundary isomorphism $\delta_{\mathcal{Q}}:K_*(\mathcal{Q}_{C^*\pi}) \cong K_{*-1}(C^*\pi)$.

(4) As $\rho_{M,N}$ is independent of the choice of $r$ and $\tilde{x}_0$, we can choose $r=0$.

By the consideration above, we get: 
$$\alpha(N)=\langle [\nu_N],[N] \rangle=\langle [\nu_N],\partial_{MV}[M]\rangle=\langle \delta_{MV}[\nu_N],[M]\rangle$$
$$2^{\epsilon}sgn(N,\nu_N)=\langle [\nu_N],2^{\epsilon}[N]_{sgn}\rangle=\langle [\nu_N],\partial_{MV}[M]_{sgn}\rangle=\langle \delta_{MV}[\nu_N],[M]_{sgn}\rangle$$

Next we have to analyse the isomorphism $\delta_{\mathcal{Q}}$, similar to Lemma 2.4 in \cite{Kubotacodim2}, we can find a representative class in $K^{-1}(N;\mathcal{Q}_{C^*\pi})$, such that it maps to $[\nu_N]$ under the isomorphism $\delta_{\mathcal{Q}}$:

\begin{Lemma}[Lemma 2.4 in \cite{Kubotacodim2}]
	\label{calkin}
	\
	
	Let $\underline{H_{C^*\pi}}$ be the trivial Hilbert $C^*\pi$ module bundle over $N$ with fiber $H_{C^*\pi}$, the standard Hilbert module. Consider the Hilbert $C^*\pi$ module bundle $\nu_N \oplus \underline{H_{C^*\pi}}$, it is trivial by Kuiper's theorem. The trivialization map $\underline{H_{C^*\pi}} \longrightarrow \nu_N \oplus \underline{H_{C^*\pi}}$ composed with the natural projection $\nu_N \oplus \underline{H_{C^*\pi}} \longrightarrow \underline{H_{C^*\pi}}$ gives us a map $\Phi:N \longrightarrow B_{C^*\pi}$. The values of the map are epimorphisms and the kernels of the map finitely generated projective. Therefore, composing with the projection map $\varPi$ to the Calkin algebra gives us a map to the invertible elements in $\mathcal{Q}_{C^*\pi}$: $\varPi \circ \Phi: N \longrightarrow GL(\mathcal{Q}_{C^*\pi})$, which represents a class 
	
	$$[\varPi \circ \Phi] \in K^{-1}(N;\mathcal{Q}_{C^*\pi})$$
	
	Then $\delta_{\mathcal{Q}}([\varPi \circ \Phi])=[\nu_N]$.
\end{Lemma}

The analysis of the map $\delta_{MV}$ also comes from Lemma 2.6 in \cite{Kubotacodim2}:
\begin{Lemma}[Lemma 2.6 in \cite{Kubotacodim2}]
	\label{MV}
	\
	
	Let $V$ be the Hilbert $C^*\pi$ module bundle over $M$ obtained by forming the trivial free $\mathcal{Q}_{C^*\pi}$ module bundles of rank one over $M_N$ and gluing along its boundary by the isomorphism $\varPi \circ \Phi$ in Lemma \ref{calkin}, that is, $V=M_N \times \mathcal{Q}_{C^*\pi}/\sim$, where $(t(x,1),T) \sim (t(x,-1),(\Pi \circ \Phi)(x)(T))$ for all $x \in N$ and $T \in \mathcal{Q}_{C^*\pi}$.
		
	Then $\delta_{MV}([\varPi \circ \Phi])=[V]$.
\end{Lemma}

Now the main theorem will follow from the following crucial theorem, which will be proved later in this section.

\begin{Theorem}
	\label{crucial}
	$V \cong \widetilde{M} \times_{\rho} \mathcal{Q}_{C^*\pi}$
\end{Theorem}

\begin{proof}[Proof of Theorem \ref{main thm}]
	To start with,
	\begin{align*}
		\rho_{M,N}\alpha(M)=\delta_{\mathcal{Q}}\rho_*\alpha(M)=\delta_{\mathcal{Q}}\rho_*\langle [\nu_M],[M] \rangle &=\langle \delta_{\mathcal{Q}}([\rho_*\nu_M]),[M] \rangle  \\
		\rho_{M,N}sgn(M',f^*\nu_M)=\delta_{\mathcal{Q}}\rho_*sgn(M',f^*\nu_M)
		&=\delta_{\mathcal{Q}}\rho_*\langle [f^*\nu_M],[M']_{sgn} \rangle \\ &=\delta_{\mathcal{Q}}\rho_*\langle [\nu_M],f_*[M']_{sgn} \rangle  \\
		&=\langle \delta_{\mathcal{Q}}([\rho_*\nu_M]),f_*[M']_{sgn} \rangle  \\
	\end{align*}
	Then,
	\begin{align*}
		\delta_{\mathcal{Q}}([\rho_*\nu_M])&=\delta_{\mathcal{Q}}([V]) \ (by \ Theorem \ \ref{crucial})\\
		&=\delta_{\mathcal{Q}}\delta_{MV}([\varPi \circ \Phi]) \ (by \ Lemma \ \ref{MV})\\
		&=\delta_{MV}([\nu_N]) \ (by \ Lemma \ \ref{calkin})\\
		\langle
		\delta_{\mathcal{Q}}([\rho_*\nu_M]),[M]
		&=\langle \delta_{MV}([\nu_N]),[M] \rangle \\
		&=\langle [\nu_N],\partial_{MV}([M]) \rangle \\
		&=\langle [\nu_N],[N] \rangle\\
		&=\alpha(N)\\
		\langle
		 \delta_{\mathcal{Q}}([\rho_*\nu_M]),f_*[M']_{sgn}
		&=\langle \delta_{MV}([\nu_N]),f_*[M']_{sgn} \rangle \\
		&=\langle [\nu_N],\partial_{MV}(f_*[M']_{sgn}) \rangle \\
		&=2^{\epsilon}\langle [\nu_N],f_*[N']_{sgn} \rangle\\
		&=2^{\epsilon}sgn(N',f^*\nu_N)\\
	\end{align*}
\end{proof}

(The exchange of $\delta_{MV}$ and $\delta_{\mathcal{Q}}$ follows from general properties of Kasparov products, see for example Theorem 2.14 in \cite{KasparovKK})

Now we begin to prove theorem \ref{crucial}, it is divided into two cases:

\subsection{Separating case}
\

Suppose $M \backslash N$ has two connected components, denote the closure of them to be $C_+$ and $C_-$ with conditions $t(N \times [0,1]) \subset C_+$ and $t(N \times [-1,0]) \subset C_-$, denote $\Gamma_+=\pi_1(C_+,x_0)$ and $\Gamma_-=\pi_1(C_-,x_0)$, then $M=C_+ \cup_N C_-$ and $\Gamma=\Gamma_+ * _{\pi} \Gamma_-$.

\begin{Lemma}
	\label{+}
	For $\forall g_+ \in \Gamma_+$, $\varphi_{g_+}(E_{[0,1]}) \subset E_{[0,1]}$ and $\varphi_{g_+}(E_{[-1,0)}) \subset E_{[-1,0)}$.
\end{Lemma}

\begin{proof}
	\
	
	Before starting the proof, let us recall some definitions in section \ref{constructiontransfer}. $E$ is the completion of $\Gamma \times_{\pi} C^*\pi=\mathop{\oplus}\limits_{s \in \Gamma/\pi}C^*\pi$. Then $\pi$ acts on $C^*\pi$ by left multiplication and induction gives a $\Gamma$ action on $\Gamma \times_{\pi} C^*\pi$. For $g \in \Gamma$, $\varphi_g \in B_{C^*\pi}$ is the continous extension of the acton of $g$ on $\Gamma \times_{\pi} C^*\pi$. For $A \subset [-1,1]$, $E_{A} \subset E$ is the subset of $E$ consisting of those elements in which $x_s=0$ for all $s \notin \chi_{\Gamma/\pi}^{-1}(A)$. Thus $\varphi_g$ is the map permuting the index $\Gamma/\pi$ by left multiplication by $g$ and we will analyze below how it changes the index.
	
	(1) Choose $s=g\pi \in \Gamma / \pi$, s.t. $\chi_{\Gamma / \pi}(s) \in (0,1]$, then $\underline{p}(g^{-1}\tilde{x}_0) \in \overline{M}_{+} \backslash N$, $g^{-1}\tilde{x}_0 \in \widetilde{M}_{+} \backslash \widetilde{N}$.
	
	Choose a loop $\gamma: [0,1] \longrightarrow C_+$ based at $x_0$ reprensenting $g_+^{-1}$ with the following properties: 
	
	$\exists \ \delta>0$, s.t. $\gamma(s)=t(x_0,s)$, $\forall s \in [0,2\delta]$ and $\gamma(s) \in C_+ \backslash N$, $\forall s \in (0,1)$
	
	Lift it to a path $\tilde{\gamma}$ in $\widetilde{M}$ with $\tilde{\gamma}(0)=g^{-1}\tilde{x}_0$, then $\tilde{\gamma}(1)=g^{-1}g_+^{-1}\tilde{x}_0$. As $\gamma((0,1)) \cap N=\emptyset$, $\tilde{\gamma}((0,1)) \cap \tilde{p}^{-1}(N)=\emptyset$, thus $\tilde{\gamma}((0,1)) \cap \widetilde{N}=\emptyset$. As $\tilde{\gamma}(0) \in \widetilde{M}_{+} \backslash \widetilde{N}$, then $\tilde{\gamma}([0,1)) \subset  \widetilde{M}_{+} \backslash \widetilde{N}$, thus $g^{-1}g_+^{-1}\tilde{x}_0=\tilde{\gamma}(1) \in \widetilde{M}_{+}$, showing that $\chi_{\Gamma / \pi}(g_+g\pi) \in [0,1]$.
	
	(2) Choose $s=g\pi \in \Gamma / \pi$, s.t. $\chi_{\Gamma / \pi}(s)=0$, then $g^{-1}\tilde{x}_0 \in \widetilde{N}$, thus $g \in \pi$.
	
	Choose a loop $\gamma$ representing $g_+^{-1}$ as in (1) and lift it to a path $\tilde{\gamma}$, then by the uniqueness of lifting, $\tilde{\gamma}(s)=\tilde{t}(g^{-1}\tilde{x}_0,s)$ for $\forall s \in [0,2\delta]$, similarly we have $\tilde{\gamma}((0,1)) \cap \widetilde{N}=\emptyset$. Then as $\tilde{\gamma}(\delta) \in \widetilde{M}_{+} \backslash \widetilde{N}$, we get $\tilde{\gamma}((0,1)) \subset  \widetilde{M}_{+} \backslash \widetilde{N}$, thus $g^{-1}g_+^{-1}\tilde{x}_0=\tilde{\gamma}(1) \in \widetilde{M}_{+}$, showing that $\chi_{\Gamma / \pi}(g_+g\pi) \in [0,1]$.
	
	(3) Choose $s=g\pi \in \Gamma / \pi$, s.t. $\chi_{\Gamma / \pi}(s) \in [-1,0)$, similarily to (1) we get $g^{-1}g_+^{-1}\tilde{x}_0 \in \widetilde{M}_{-}$ and  $\chi_{\Gamma / \pi}(g_+g \pi) \in [-1,0]$. If $\chi_{\Gamma / \pi}(g_+g\pi)=0$, then $g_+g \in \pi \subset \Gamma_+$, but $g_+ \in \Gamma_+$, we get $g \in \Gamma_+$. Then by argument in (2) $g^{-1}\tilde{x}_0 \in \widetilde{M}_{+}$, contradiction with $\chi_{\Gamma / \pi}(s) \in [-1,0)$. Thus $\chi_{\Gamma / \pi}(g_+g\pi) \in [-1,0)$.
	
	Combining the arguments in the three cases, we get the result.
	
\end{proof}

Similar arguments shows that:

\begin{Lemma}
	\label{-}
	For $\forall g_- \in \Gamma_-$, $\varphi_{g_-}(E_{[-1,0]}) \subset E_{[-1,0]}$ and $\varphi_{g_-}(E_{(0,1]}) \subset E_{(0,1]}$.
\end{Lemma}

Combining the two lemmas we get:

\begin{Lemma}
	\label{0}
	For $\forall g \in \pi$, $\varphi_{g}(E_{(-1,0]}) \subset E_{(-1,0]}$, $\varphi_{g}(E_{(0,1]}) \subset E_{(0,1]}$ and $\varphi_{g}(E_{\{0\}}) \subset E_{\{0\}}$. Furthermore, $\varphi_g|_{E_{\{0\}}}=\lambda_{g}: E_{\{0\}}=C^*\pi \longrightarrow E_{\{0\}}=C^*\pi$, where $\lambda_{g}$ is the left multiplication of $g$ on $C^*\pi$.
\end{Lemma}

Before starting to prove theorem \ref{crucial}, let us recall some constructions in \cite{Kubotacodim2}. Let $H$ be a Hilbert $C^*\pi$ module bundle, whose typical fiber is the standard Hilbert module $H_{C^*\pi}$ and with structure group $GL_{C^*\pi}(H_{C^*\pi})$ with norm topology. We write $GL(H)$  for corresponding right principal $GL_{C^*\pi}(H_{C^*\pi})$–bundle with fiber $GL_{C^*\pi}(H_{C^*\pi},H_x)$. The products $(A,T) \mapsto AT$ and $(A,[T]) \mapsto [AT]$ induce left $GL_{C^*\pi}(H_{C^*\pi})$ actions on $B_{C^*\pi}$ and $\mathcal{Q}_{C^*\pi}$. We have the associated bundles $B(H)=GL(H) \times_{GL_{C^*\pi}(H_{C^*\pi})} B_{C^*\pi}$ and $\mathcal{Q}(H)=GL(H) \times_{GL_{C^*\pi}(H_{C^*\pi})} \mathcal{Q}_{C^*\pi}$. Then the following are easily verified: (See Remark 2.16 in \cite{Kubotacodim2}.)

(1) Let $X$ be a path connected topological space and $\widetilde{X}$ be its universal cover. Let $\rho:\pi_1(X) \longrightarrow GL_{C^*\pi}(H_{C^*\pi})$ be a representation, then the associated bundle $\widetilde{X} \times_{\pi_1(X)} B_{C^*\pi}$ is isomorphic to $B(\widetilde{X} \times_{\pi_1(X)} H_{C^*\pi})$.

(2) A bounded bundle map $T: H \longrightarrow H'$ induces bundle maps $B(H) \longrightarrow B(H')$ and $\mathcal{Q}(H) \longrightarrow \mathcal{Q}(H')$.

(3) A bounded bundle map $T: X \times H_{C^*\pi} \longrightarrow X \times H_{C^*\pi}$ between trival bundles is identified with a continuous function $X \longrightarrow B_{C^*\pi}$. The bundle maps on $B(X \times H_{C^*\pi})=X \times B_{C^*\pi}$ and $\mathcal{Q}(X \times H_{C^*\pi})=X \times \mathcal{Q}_{C^*\pi}$ induced from T are the multiplication of $T$ from the left.

Now we start the proof. Let $U=\widetilde{M} \times_{\rho} \mathcal{Q}_{C^*\pi}$, because of lemma \ref{+} $\varphi$ defines a group homomorphism $\varphi_+:\Gamma_+ \longrightarrow GL_{C^*\pi}(E_{[0,1]})$. By definition $\varPi(\varphi_+ \oplus Id)=\rho|_{\Gamma_+}$. Thus if we set $\widetilde{C_+}$ to be the universal cover of $C_+$ and $W_+=\widetilde{C_+} \times_{\varphi_+ \oplus Id} (E_{[0,1]} \oplus H_{C^*\pi})$, a Hilbert $C^*\pi$ module bundle over $C_+$, then $U|_{C_+} \cong \mathcal{Q}(W_+)$.

By lemma \ref{0}, $\varphi_+|_{\pi}=\lambda_{\pi} \oplus \varphi_0:\pi \longrightarrow GL_{C^*\pi}(E_{\{0\}} \oplus E_{(0,1]})$ for some $\varphi_0$. If we set $W_0=\widetilde{N} \times_{\varphi_0 \oplus Id} (E_{(0,1]} \oplus H_{C^*\pi})$, then $W_+|_N \cong \nu_N \oplus W_0$.

By lemma \ref{-}, $\varphi$ defines a group homomorphism $\varphi_-: \Gamma_- \longrightarrow GL_{C^*\pi}(E_{(0,1]})$ by restriction. If we set $\widetilde{C_-}$ to be the universal cover of $C_-$ and $W_-=\widetilde{C_-} \times_{\varphi_- \oplus Id} (E_{(0,1]} \oplus H_{C^*\pi})$, a Hilbert $C^*\pi$ module bundle over $C_-$, then we can glue $\mathcal{Q}(W_+)$ and $\mathcal{Q}(W_-)$ by the bundle isomorphism $\mathcal{Q}(p_1):\mathcal{Q}(W_+)|_N \longrightarrow \mathcal{Q}(W_-)|_N$ induced by projection:
$$p_1: W_+|_N \cong \nu_N \oplus W_0 \longrightarrow W_0=W_-|_N$$

Denote $\iota:E_{(0,1]} \oplus H_{C^*\pi} \longrightarrow E_{[0,1]} \oplus H_{C^*\pi}$ and $\mathfrak{p}:E_{[0,1]} \oplus H_{C^*\pi} \longrightarrow E_{(0,1]} \oplus H_{C^*\pi}$ to be the inclusion and projection map respectively. Then $\mathfrak{p}_*:\mathcal{Q}_{C^*\pi}(E_{[0,1]} \oplus H_{C^*\pi}) \longrightarrow \mathcal{Q}_{C^*\pi}(E_{(0,1]} \oplus H_{C^*\pi}), \ \mathfrak{p}_*([T])=[\mathfrak{p}T\iota]$ is an isomorphism with the inverse $\mathfrak{p}_*^{-1}([T'])=[\iota T' \mathfrak{p}]$. Since $W_+$ and $W_-$ are flat bundles, $\mathcal{Q}(W_+)|_N$ and $\mathcal{Q}(W_-)|_N$ are also flat bundles, thus their pull back to $\widetilde{N}$ are trivial. Note that we have the following commutative diagram: 

$$
\begin{tikzcd}
	\widetilde{N} \times \mathcal{Q}_{C^*\pi}(E_{[0,1]} \oplus H_{C^*\pi}) \rar{(id,\mathfrak{p}_*)}[swap]{\cong} \dar & \widetilde{N} \times \mathcal{Q}_{C^*\pi}(E_{(0,1]} \oplus H_{C^*\pi}) \dar \\
	\mathcal{Q}(W_+) \rar{\mathcal{Q}(p_1)}[swap]{\cong} & \mathcal{Q}(W_-) \\
\end{tikzcd}
$$

Where the vertical maps are the natural bundle maps for the pull back $\tilde{p}:\widetilde{N} \longrightarrow N$.
Thus the glued bundle $Z$ is also flat. Suppose the honolomy of $Z$ is $\rho_Z$, we claim that $\rho_Z=\rho$. To prove the claim, note first that by the gluing construction, $\rho_Z(g_+)=\varPi(\varphi_+(g_+) \oplus Id)=\rho_{g_+}$ for all $g_+ \in \Gamma_+$ and $\rho_Z(g_-)=\mathfrak{p}_*^{-1}(\varPi(\varphi_-(g_-) \oplus Id))$ for all $g_- \in \Gamma_-$. Since $\Gamma=\Gamma_+*_{\pi}\Gamma_-$, it suffices to prove $\mathfrak{p}_*^{-1}(\varPi(\varphi_-(g_-) \oplus Id))=\rho_{g_-}$ for all $g_- \in \Gamma_-$. Note that $\mathfrak{p}_*^{-1}(\varPi(\varphi_-(g_-) \oplus Id))=[\iota (\varphi_-(g_-) \oplus Id)\mathfrak{p}]$. For any $z \in E_{[0,1]} \oplus H_{C^*\pi}=E_0 \oplus E_{(0,1]} \oplus H_{C^*\pi}$, write $z=(z_1,z_2,z_3)$, then

\begin{align*}
	\iota (\varphi_-(g_-) \oplus Id)\mathfrak{p}(z)&=\iota (\varphi_-(g_-) \oplus Id)(z_2,z_3) \\ &=\iota (\varphi_-(g_-)(z_2),z_3) \\ &=\iota (\varphi_{g_-}(z_2),z_3) \ (\text{By Definition of } \varphi_-) \\ &=(0,\varphi_{g_-}(z_2),z_3)
\end{align*}
\begin{align*}
	(P_{[0,1]} \varphi_{g_-}  I_{[0,1]} \oplus Id)(z)&=(P_{[0,1]} \varphi_{g_-}  I_{[0,1]} \oplus Id)(z_1,z_2,z_3)\\ &=(P_{[0,1]} \varphi_{g_-}  I_{[0,1]} (z_1,z_2),z_3) \\ &=(P_{[0,1]} \varphi_{g_-}  I_{[0,1]} (z_1,0),0)+(P_{[0,1]} \varphi_{g_-}  I_{[0,1]} (0,z_2),z_3)\\ &=(P_{[0,1]} \varphi_{g_-}  I_{[0,1]} (z_1,0),0)+(0,\varphi_{g_-}(z_2),z_3) \\ & \ \text{By Lemma } \ref{-}
\end{align*}
 
 Since $E$ is the completion of $\mathop{\oplus}\limits_{s \in \Gamma/\pi}C^*\pi$ and $\varphi_{g_-}$ is given by permuting the index $\Gamma/\pi$, the difference of $\iota (\varphi_-(g_-) \oplus Id)\mathfrak{p}$ and $P_{[0,1]} \varphi_{g_-}  I_{[0,1]} \oplus Id$ is a rank one operator, then their image in the Calkin algebra is the same. Recall that $\rho_{g_-}$ is the image of $P_{[0,1]} \varphi_{g_-}  I_{[0,1]} \oplus Id$ in the Calkin algebra, thus $\rho_Z(g_-)=\rho_{g_-}$ for all $g_- \in \Gamma_-$. We have proven the claim, then $Z$ is isomorphic to $U$ as it has the same holonomy as $U$.

Now by Kuiper's Theorem, $W_+$ and $W_-$ are trivial bundles with trivialization:
$$v_+: W_+ \longrightarrow C_+ \times H_{C^*\pi} \ \ \  v_-: W_- \longrightarrow C_- \times H_{C^*\pi}$$

Then the map $v_-|_N \circ p_1 \circ v_+^{-1}|_N(x,a)=(x,\Phi_x(a))$ defines a map $\Phi: N \longrightarrow B_{C^*\pi}$, which can be chosen as the map $\Phi$ in lemma \ref{calkin}. Thus by the construction of $V$ in lemma \ref{MV}, $V=Z \cong U=\widetilde{M} \times_{\rho} \mathcal{Q}_{C^*\pi}$.

\subsection{Non-separating case}
\

Suppose $M \backslash N$ is connected, $\partial M_N$ is the disjoint union of two copies of $N$ and we denote $N_+=t(N \times \{1\})$ and $N_-=t(N \times \{-1\})$, then $\partial M_N=N_+ \coprod N_-$. Set $x_1=t(x_0,1) \in N_+$ and $x_{-1}=t(x_0,-1)$. Choose any path $\gamma_0$ in $M_N$ from $x_{-1}$ to $x_{1}$ and join it with the path $\gamma_1:[0,1] \longrightarrow M, \gamma_1(s)=t(x_0,2s-1), \text{for } s \in [0,1]$, we get a loop $\gamma$ based at $x_0$. We define $a \in \Gamma$ to be the element represented by $\gamma$ and denote $\pi_1(M_N,x_1)=G$. We define $i_+:G \longrightarrow \Gamma$ to be the map induced by the inclusion map and $i_-:G \longrightarrow \Gamma$ to be the composition of the isomorphism $\gamma_0^{\sharp}:\pi_1(M_N,x_{-1}) \cong \pi_1(M_N,x_1)=G$ induced by $\gamma_0$ and the induced map of the inclusion map $M_N \longrightarrow M$, then $\Gamma=G*\mathbb{Z}/<<ai_+(g)a^{-1}i_-(g^{-1})>>$, where $<<ai_+(g)a^{-1}i_-(g^{-1})>>$ is the normal closure of the subset $\{ai_+(g)ai_-(g^{-1})\}$. Moreover, $i_+$ is the  inclusion map under this identification. If we set $r:M_N \longrightarrow M$ to be the map defined by \begin{equation*}
	r(y)=\left\{
	\begin{aligned}
		& t(x,2s+2)  &y&=t(x,s) \in t(N \times [-2,-1])\\
		& t(x,2s-2)  &y&=t(x,s) \in t(N \times [1,2])\\
		& y &y& \in M_N \backslash t(N \times [-2,2])
	\end{aligned}
	\right.
\end{equation*}

Then $r_*:G \longrightarrow \Gamma$ is the same as $i_+$.

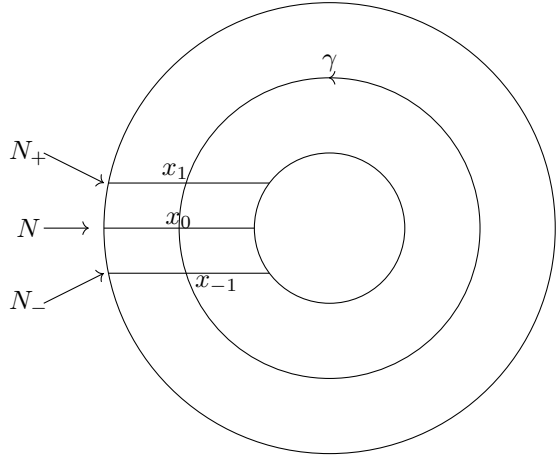
$$
\begin{tikzpicture}
	\draw (0,0) circle (3);
	\draw (0,0) circle (2);
	\draw (0,0) circle (1);
	\draw (-3,0)--(-1,0);
	\draw (-2.94,0.6)--(-0.8,0.6);
	\draw (-2.94,-0.6)--(-0.8,-0.6);
	\draw[->] (-3.8,0)--(-3.2,0);
	\draw[->] (-3.8,1)--(-3,0.6);
	\draw[->] (-3.8,-1)--(-3,-0.6);
	\draw[->] (0.01,2)--(-0.01,2);
	\node at (-4,0) {$N$};
	\node at (0,2.2) {$\gamma$};
	\node at (-4,1) {$N_+$};
	\node at (-4,-1) {$N_-$};
	\node at (-2,0.1) {$x_0$};
	\node at (-2.05,0.75) {$x_1$};
	\node at (-1.5,-0.75) {$x_{-1}$};
\end{tikzpicture}
$$

Similar to the proofs of Lemma \ref{+}, we can prove the following lemmas:

\begin{Lemma}
	For $\forall g \in G$, $\varphi_g(E_{[0,1]}) \subset E_{[0,1]}$.
\end{Lemma}

\begin{Lemma}
	\label{2+}
	
	For $\forall h \in \pi$, $\varphi_{i_+h}(E_{0}) \subset E_{0}$, $\varphi_{i_+h}(E_{(0,1]}) \subset E_{(0,1]}$, and $\varphi_{i_+h}|_{E_{0}}=\lambda_{h}$.
	
\end{Lemma}

\begin{Lemma}
	
	\label{t}
	
	$\varphi_a(E_{(0,1]}) \subset E_{[0,1]}$, $\varphi_a(E_0) \subset E_{[-1,0)}$, $\varphi_{a^{-1}}(E_{[0,1]}) \subset E_{(0,1]}$.
\end{Lemma}

Set $U=\widetilde{M} \times_{\rho} \mathcal{Q}_{C^*\pi}$ and $W=\widetilde{M_N} \times_{\varphi \circ i_+ \oplus Id} (E_{[0,1]} \oplus H_{C^*\pi})$. Since $r_*=i_+$ and by definition $\varPi(\varphi \circ i_+ \oplus Id)=\rho|_{G}$, we get $r^*U=\mathcal{Q}(W)$.

Set $W_0=\widetilde{N} \times_{\varphi \circ i_{+} \oplus Id} (E_{(0,1]} \oplus H_{C^*\pi})$, note that $W|_{N_{-}}=\widetilde{N} \times_{\varphi \circ i_{-} \oplus Id} (E_{[0,1]} \oplus H_{C^*\pi})$ and by lemma \ref{2+}, $W|_{N_{+}}=\widetilde{N} \times_{\varphi \circ i_{+} \oplus Id} (E_{[0,1]} \oplus H_{C^*\pi}) \cong \nu_N \oplus W_0$.

By lemma \ref{t}, $\varphi_{a^{-1}}:E_{[0,1]} \longrightarrow E_{(0,1]}$ is an isomorphism. As $\varphi$ is a homomorphism, $\varphi_{i_-h}=\varphi_a\varphi_{i_+h}\varphi_{a^{-1}}$ for $\forall h \in \pi$, thus $W|_{N_{-}} \cong W_0$. If we identify $\mathcal{Q}(W)|_{N_+}$ with $\mathcal{Q}(W)|_{N_-}$ by the bundle map induced by projection:
$$p_1: W|_{N_{+}} \cong \nu_N \oplus W_0 \longrightarrow W_0 \cong W|_{N_{-}}$$

then the glued bundle $Z$ is isomorphic to $U$ as it has the same holonomy as $U$.

Now by Kuiper's Theorem, $W$ is a trivial bundle with trivialization:
$$v: W \longrightarrow M_N \times H_{C^*\pi}$$

Then the map $v|_{N_-} \circ p_1 \circ v^{-1}|_{N_+}(x,a)=(x,\Phi_x(a))$ defines a map $\Phi: N \longrightarrow B_{C^*\pi}$, which can be chosen as the map $\Phi$ in lemma \ref{calkin}. Thus by the construction of $V$ in lemma \ref{MV}, $V=Z \cong U=\widetilde{M} \times_{\rho} \mathcal{Q}_{C^*\pi}$.

	\section{Proof of Lemma \ref{indelem}}
\label{pfindelem}
In this section, we will give the proof of lemma \ref{indelem}, the independence of choices in the definition of the transfer map.

To begin with, for any $r \in (-1,1)$, let $S_r=\{\gamma\pi| \ \chi(\underline{p}(\gamma^{-1}\tilde{x}_0)) \geq r \} \subset \Gamma/\pi$. The Hilbert module used in the definition of $\rho_{M,N}$ is $E_{S_r}$, the following lemma analyses difference of the sets under different choices:

\begin{Lemma}
	\label{diff}
	For any $r,r' \in (-1,1)$, $S_{r} \bigtriangleup S_{r'}$ is a finite set.
\end{Lemma}

\begin{proof}
	Without loss of generosity we can assume $r \leq r'$, 
	
	$
	\begin{aligned}
		S_{r} \bigtriangleup S_{r'}=\{\gamma\pi| \ r \leq \chi(\underline{p}(\gamma^{-1}\tilde{x}_0)< r'\}&=\{\gamma\pi| \ \underline{p}(\gamma^{-1}\tilde{x}_0) \in \bar{t}(N \times [r,r')) \} \\
		& \subset \{\gamma\pi| \  \underline{p}(\gamma^{-1}\tilde{x}_0) \in \bar{t}(N \times [r,r']) \}
	\end{aligned}
	$
	
	As $\underline{p}(\gamma^{-1}\tilde{x}_0)$ is discrete and $\bar{t}(N \times [r,r'])$ is compact, $S_{r} \bigtriangleup S_{r'}$ is a finite set.
\end{proof}

Now for any $r,r' \in (-1,1)$, suppose $\rho,\rho': C^*\Gamma \longrightarrow B_{C^*\pi}$ are the maps constructed in section 3 corresponding to the data $r$ and $r'$ respectively. Denote $S=S_{r}$, $S'=S_{r'}$ $S_0=S \cap S'$, $S_1=S \backslash S'$, $S_1'=S' \backslash S$, then $S_1$ and $S_1'$ are finite sets. Suppose $|S_1|=k, \ |S_1'|=k'$.

\begin{proof}[Proof of Lemma \ref{indelem}]
	\
	
	Fix an isomorphism $f:E_{S_0} \oplus H_{C^*\pi} \cong H_{C^*{\pi}}$, then we can define isomorphisms:
	
	$q:E_{S}=E_{S_1} \oplus E_{S_0} \oplus H_{C^*\pi} \stackrel{Id \oplus f}{\cong} (C^*\pi)^k \oplus H_{C^*\pi} \stackrel{h}{\cong} H_{C^*\pi}$, 
	
	where $h:(x_1,x_2,...,x_k),(y_1,y_2,....) \mapsto (x_1,x_2,...,x_k,y_1,...)$.

	$q':E_{S'}=E_{S_1'} \oplus E_{S_0} \oplus H_{C^*\pi} \stackrel{Id \oplus f}{\cong} (C^*\pi)^{k'}\oplus H_{C^*\pi} \stackrel{h'}{\cong} H_{C^*\pi}$, 
	
	where $h':(x_1,x_2,...,x_{k'}),(y_1,y_2,....) \mapsto (x_1,x_2,...,x_{k'},y_1,...)$.

	For any $g \in G$, recall that $\rho_g$ is the element $\phi_g=P_{S}\varphi_gI_{S} \oplus Id \in B_{C^*\pi}(E_{S} \oplus H_{C^*\pi}) \stackrel{q_*}{\cong} B_{C^*\pi}$ projected onto $\mathcal{Q}_{C^*\pi}$, $\rho_g' $ is the element $\phi'_g=P_{S'}\varphi_gI_{S'} \oplus Id \in B_{C^*\pi}(E_{S'} \oplus H_{C^*\pi}) \stackrel{q'_*}{\cong} B_{C^*\pi}$ projected onto $\mathcal{Q}_{C^*\pi}$. Denote $\psi_g=P_{S_0}\varphi_gI_{S_0} \oplus 0 \in B_{C^*\pi}(E_{S_0} \oplus H_{C^*\pi}) \stackrel{f_*}{\cong} B_{C^*\pi}$. In the following
	we determine the relation of $\psi_g$ to $\phi_g$ and $\phi'_g$.
	
	Consider the composition $E_S \oplus H_{C^*\pi}=E_{S_1} \oplus (E_{S_0} \oplus H_{C^*\pi})$, we claim that there is a homomorphism $s:E_{S_0} \oplus H_{C^*\pi} \longrightarrow E_{S_1}$, such that $\phi_g(0,x)=(s(x),\psi_g(x))$ for all $x \in E_{S_0} \oplus H_{C^*\pi}$. To prove the claim, we write $x=(y,z) \in E_{S_0} \oplus H_{C^*\pi}$, then,
	
	\begin{align*}
		\phi_g(0,x)=(P_S\varphi_gI_S(0,y),z)&=(P_S\varphi_gI_{S_0}(y),z) \\ &=(P_{S_1}\varphi_gI_{S_0}(y),P_{S_0}\varphi_gI_{S_0}(y),z) \\ &\in E_{S_1} \oplus E_{S_0} \oplus H_{C^*\pi}
	\end{align*}
	
 	Note that $\psi_g(x)=(P_{S_0}\varphi_gI_{S_0}(y),z) \in  E_{S_0} \oplus H_{C^*\pi}$, thus if we set $s(x)=P_{S_1}\varphi_gI_{S_0}(y)$ for $x=(y,z) \in E_{S_0} \oplus H_{C^*\pi}$, we will get the equation in the claim.
 	
 	Next we consider the relation of $f_*\psi_g$ to $q_*\phi_g$ and $q_*\phi'_g$ in $B_{C^*\pi}$. Denote $R,L$ to be the right and left shift operator on $H_{C^*\pi}$, then for any $y \in H_{C^*\pi}$, $h(0,y)=R^ky, L^kh(z,y)=y, \forall z \in (C^*\pi)^k$. For any $y \in H_{C^*\pi}$,
 	
 	\begin{align*}
 		q_*\phi_g(R^ky)=q\phi_gq^{-1}(R^ky)&=h(Id \oplus f)\phi_g(Id \oplus f^{-1})h^{-1}(R^ky) \\ &=h(Id \oplus f)\phi_g(0,f^{-1}(y)) \\ &=h(Id \oplus f)(sf^{-1}(y),\psi_gf^{-1}(y)) \\ &=h(sf^{-1}(y),f_*\psi_g(y)) 
 	\end{align*}
 	
 	Thus $L^kq_*\phi_g(R^ky)=f_*\psi_g(y)$, i.e. $ L^k(q_*\phi_g)R^k=f_*\psi_g$. Similarly we get $L^{k'}(q'_*\phi'_g)R^{k'}=f_*\psi_g$.
	
	As $L^*=R,LR=Id,RL=Id-P_1$, where $P_1$ is the projection onto the first component. Since $L,R$ are unitaries in $\mathcal{Q}_{C^*\pi}$ with $L^{-1}=R$ and $\rho_g=\varPi q_*\phi_g$, $\rho'_g=\varPi q'_*\phi'_g$, by the arguments above $\rho$ and $\rho'$ are unitary equivalent, then $\rho_*=\rho'_*$. We finished the proof of Lemma \ref{indelem}.

\end{proof}

	\section{Extension to reduced $C^*$ algebra}
In this section, we will improve the transfer map constructed of section \ref{transfer map} to the case of reduced $C^*$-algebras. We write $C^*_rG$ for the reduced group $C^*$-algebra of a group $G$.

We can define a group homomorphism $\rho: \Gamma \longrightarrow \mathcal{Q}_{C^*_r\pi}$ in the same way as in section \ref{transfer map}, i.e. $\rho_g=\varPi \circ (P_{[r,1]}\varphi_gI_{[r,1]} \oplus 0)$ with $\rho_g^*=\rho_{g^{-1}}$. The extension of the map to a homomorphism from $C^*_r\Gamma$ to $\mathcal{Q}_{C^*_r\pi}$ is guaranteed by the following control lemma:

\begin{Lemma}
	For any $\sum\limits_{i=1}^k a_ig_i \in \mathbb{C}\Gamma$, $\lVert\sum\limits_{i=1}^k a_i\varphi_{g_i}\rVert_{B_{C^*_r\pi}(E)} \leq \lVert\sum\limits_{i=1}^k a_ig_i \rVert_{C^*_r\Gamma}$
\end{Lemma}

\begin{proof}
	For a discrete group $G$, recall that $C^*_rG$ is defined to be the closure of $\mathbb{C}G$ in bounded operators on $l^2(G)$. Denote $m_g$ to the operator of the action of $g \in G$ on $l^2(G)$. Using a transversal of $\pi$ in $\Gamma$, identity $l^2(\Gamma)$
	with $l^2(\Gamma/\pi,l^2(\pi))$. Given any element $v \in l^2(\pi)$, obsevere that there is the following map between $E$ and $l^2(\Gamma/\pi,l^2(\pi))$ given by:
	
	$x=(x_s)_{s \in \Gamma/\pi} \in E \mapsto x(v) \stackrel{def}{=} (x_s(v))_{s \in \Gamma/\pi} \in l^2(\Gamma/\pi,l^2(\pi))$
	
	The map has the following three properties:
	
	(1) $\lVert x(v) \rVert \leq \lVert x \rVert_E \lVert v \rVert_{l^2(\pi)}$, in particular, the map is continous.
	
	Proof: 
	\begin{align*}
		\lVert x(v) \rVert ^2= \sum\limits_{s \in \Gamma/\pi} \lVert x_s(v) \rVert_{l^2(\pi)} ^2&=\sum\limits_{s \in \Gamma/\pi} (x_s(v),x_s(v))_{l^2(\pi)} \\ &=\sum\limits_{s \in \Gamma/\pi} (x_s^*x_s(v),v)_{l^2(\pi)} \\ & \leq \lVert \sum\limits_{s \in \Gamma/\pi}x_s^*x_s \rVert_{B_{l^2(\pi)}} \lVert v \rVert_{l^2(\pi)} ^2 \\ &=\lVert x \rVert_E^2  \lVert v \rVert_{l^2(\pi)}^2 
	\end{align*}

	Thus $\lVert x(v) \rVert \leq \lVert x \rVert_E \lVert v \rVert_{l^2(\pi)}$.
	
	(2) The map is $\Gamma$ equivariant: $(\varphi_gx)(v)=m_g(x(v))$ for all $g \in \Gamma$.
	
	(3) $\lVert x \rVert_E=\sup\limits_{\substack{v \in l^2(\pi) \\ \lVert v \rVert=1}}\lVert x(v) \rVert$.
	
	 Proof: Observe that $\sum\limits_{s \in \Gamma/\pi}x_s^*x_s$ is self-adjoint, thus
	 \begin{align*}
	 	\lVert \sum\limits_{s \in \Gamma/\pi}x_s^*x_s \rVert=\sup\limits_{\substack{v \in l^2(\pi) \\ \lVert v \rVert=1}}(\sum\limits_{s \in \Gamma/\pi}x_s^*x_sv,v)_{l^2(\pi)}
	 \end{align*}
 
 	 Thus $\lVert x \rVert_E^2=\lVert \sum\limits_{s \in \Gamma/\pi}x_s^*x_s \rVert=\sup\limits_{\substack{v \in l^2(\pi) \\ \lVert v \rVert=1}}(\sum\limits_{s \in \Gamma/\pi}x_s^*x_sv,v)_{l^2(\pi)}=\sup\limits_{\substack{v \in l^2(\pi) \\ \lVert v \rVert=1}}\lVert x(v) \rVert^2$, giving the result.
	
	To prove the lemma, we note that by the properties above, we have 
	\begin{align*}
		\lVert \sum\limits_{i=1}^k a_i\varphi_{g_i} \rVert_{B_{C^*_r\pi}(E)}=\sup\limits_{\substack{x \in E \\ \lVert x \rVert=1}} \lVert \sum\limits_{i=1}^k a_i\varphi_{g_i}x \rVert_E &\stackrel{(3)}{=}\sup\limits_{\substack{x \in E \\ \lVert x \rVert=1}} \sup\limits_{\substack{v \in l^2(\pi) \\ \lVert v \rVert=1}} \lVert \sum\limits_{i=1}^k a_i\varphi_{g_i}x(v) \rVert \\ &\stackrel{(2)}{=}\sup\limits_{\substack{x \in E \\ \lVert x \rVert=1}} \sup\limits_{\substack{v \in l^2(\pi) \\ \lVert v \rVert=1}} \lVert \sum\limits_{i=1}^k a_im_{g_i}(x(v)) \rVert \\
	\end{align*}

	As $\lVert \sum\limits_{i=1}^k a_im_{g_i}(x(v)) \rVert \leq \lVert \sum\limits_{i=1}^k a_im_{g_i} \rVert \lVert x(v) \rVert \stackrel{(1)}{\leq} \lVert \sum\limits_{i=1}^k a_ig_i \rVert_{C^*_r\Gamma} \lVert x \rVert_E \lVert v \rVert_{l^2(\pi)}$, we get $\lVert\sum\limits_{i=1}^k a_i\varphi_{g_i}\rVert_{B_{C^*_r\pi}(E)} \leq \lVert\sum\limits_{i=1}^k a_ig_i \rVert_{C^*_r\Gamma}$.
\end{proof}

Now by the lemma above we can extend $\varphi$ to a continous homomorphism $\varphi:C^*_r \Gamma \longrightarrow B_{C^*\pi}(E)$, then $\rho=\varPi \circ (P_{[r,1]}\varphi I_{[r,1]} \oplus 0): C^*_r\Gamma \longrightarrow \mathcal{Q}_{C^*_r\pi}$ gives a definition of transfer map in reduced $C^*$-algebra case.

	\section{Further directions}
We will talk about possible further directions in this section. We have constructed a transfer map in section \ref{transfer map}. Using the desciption of KK theory by extension of $C^*$-algebras, the construction is connected to the KK theoretic construction of transfer maps in \cite{generaltransfer}. Thus it provides a direct and explict way of understanding the K theoretic transfer maps in the codimension 1 case.

Moreover, there is a way to construct L theoretic signatures by using the chain complexes. The signature of a closed manifold $P$ lies in an abelian group $L^n(\mathbb{Z}\pi_1(P))$, known as the symmetric L group. By using suspension rings, one can establish analogusly a transfer map for symmetric L theory, thereby establishing a relationship between the L-theoretic indexes of the corresponding manifolds. Additionally, there is also another kind of L group used in surgery theory, called the quadratic L group. The two L groups are related by the symmetrilization map $1+T:L_n(\mathbb{Z}\pi_1(P)) \longrightarrow L^n(\mathbb{Z}\pi_1(P))$. The construction can be extended similarly to quadratic L theory and such construction is related to the long exact sequence of quadratic L theory of codimension 1 submanifolds. (For the long exact sequnence see case (A) and (B) in Chapter 12 of \cite{Wallsurgery})

It is an interesting question to clarify all the relationships of the codimension 1 transfer maps constructed in K and L theory.
	\endgroup
	\bibliographystyle{abbrv}
	\bibliography{refer}
\end{document}